\documentclass{article}

\usepackage{amssymb, amsmath}
\usepackage{cite}

\begin{document}


\begin{center}
   \LARGE \textbf{Schr\"{o}dinger-improved Boussinesq system in two space dimensions}
\end{center}
\vspace{5pt}
\begin{center}
   \large Tohru Ozawa and Kenta Tomioka \\
Department of Applied Physics, Waseda University   \\
Tokyo 169-8555, Japan
\end{center}
\vspace{5pt}
\begin{center}
   \textbf{Abstract.}
\end{center}

We study the Cauchy problem for the Schr\"{o}dinger-improved Boussinesq system in a two dimentional domain.
Under natural assumptions on the data without smallness, we prove the existence and uniqueness of global strong solutions.
Moreover, we consider the vanishing ''improvement'' limit of global solusions as the coefficient of the linear term of the highest order in the equation of ion sound waves tends to zero.
Under the same smallness assumption on the data as in the Zakharov case, solutions in the vanishing ''improvement'' limit are shown to satisfy the Zakharov system.  \\



\pagebreak


\section{Introduction.}

\ \ \ In this paper we study the Cauchy problem for the Schr\"{o}dinger-improved Boussinesq system (S-iB)
\begin{align*}
   \text{(S-iB)}
   \begin{cases}
   i \partial_t u + \Delta u = vu,    \\
   \partial_t^2 v - \Delta v - \Delta \partial_t^2 v = \Delta |u|^2   
   \end{cases}
\end{align*}
with initial data
\begin{align*}
   (u(0),v(0),\partial_t v(0)) = (\varphi,\psi_0,\psi_1)   
\end{align*}
given at $t=0$, where $u:\mathbb R \times \Omega \to \mathbb C, v:\mathbb R \times \Omega \to \mathbb R$, and $\Omega \subset \mathbb R^2$ is a domain with smooth boundary $\partial \Omega$.
The Laplacian $\Delta$ is understood to be the self-adjoint realization in the Hilbert space $L^2(\Omega)$ with domain $D(\Delta)=(H^2\cap H_0^1)(\Omega)$, where $H^2(\Omega)$ is the Sobolev space of the second order and $H_0^1(\Omega)$ is the Sobolev space of the first order with vanishing condition on $\partial \Omega$.
The system (S-iB) is regarded as a substitute for the Zakharov system (Z)
\begin{align*}
   \text{(Z)}
   \begin{cases}
   i \partial_t u + \Delta u = vu,    \\
   \partial_t^2 v - \Delta v = \Delta |u|^2   
   \end{cases}
\end{align*}
and the ''improvement'' has been made on the second equation of (S-iB) by modifying the dispersion relation of ion-sound waves of Boussinesq type.
See \cite{M} for details.

In spite of a large literature on (S-iB) and (Z) in the whole space-time $\mathbb R \times \mathbb R^n$ (see \cite{A,CO2,EP,FO,OT2,WC1,WC2,WC3,WC4} for (S-iB) and \cite{AA1,AA2,BHHT,BC,FPZ,GTV,GGKZ,GLNW,GN,GZP,KT,
KPV,K,MN1,MN2,OT3,OT4,SW,SS1,SS2,T} for (Z)), there are few papers treating those systems in $\mathbb R \times \Omega$ (see \cite{AA1,OT1,YGN,YN1,YN2} for (Z)).
A major reason consists in the lack of the Fourier transform, by which Strichartz estimates, Bourgain's method, and I-method are available.

The first purpose of this paper is to prove the existence and uniqueness of global strong solutions to (S-iB) in $\mathbb R \times \Omega$.
We prove:   \\

\noindent
\textbf{Theorem 1.}\ \ \textit{Let $(\varphi,\psi_0,\psi_1) \in D(\Delta) \oplus H_0^1(\Omega) \oplus (H_0^1 \cap (-\Delta)^{1/2}L^2)(\Omega)$. 
Then}:   \\
\noindent
(1) (S-iB) \textit{has a unique solution $(u,v,\partial_t v)$ with} 
\begin{align*}
   &u \in (L_{loc}^{\infty} \cap C_w) (\mathbb R ; D(\Delta)) \cap C(\mathbb R ; H_0^1(\Omega)) \cap C^1(\mathbb R ; H^{-1}(\Omega)),   \\
   &\partial_t u \in (L_{loc}^{\infty} \cap C_w) (\mathbb R ; L^2(\Omega) )   \\
   &v \in C^2(\mathbb R ; H_0^1(\Omega)),   \\
   &(-\Delta)^{-1/2} \partial_t v \in C(\mathbb R ; L^2(\Omega)),
\end{align*}
\textit{and $(u(0), v(0), \partial_t v(0))=(\varphi,\psi_0,\psi_1)$, where $H^{-1}(\Omega)=(H_0^1(\Omega))^*$.}   \\
(2) \textit{The following total charge and energy are conserved in time} :
\begin{align*}
   &\|u(t)\|_2^2,   \\
   &E(t) := \|\nabla u(t)\|_2^2 + \frac{1}{2} (\|v(t)\|_2^2 
   + \|\partial_t v(t)\|_2^2 + \|(-\Delta)^{-1/2}\partial_t v(t)\|_2^2) + (v(t)||u(t)|^2),
\end{align*}
\textit{where $\|\cdot\|_2$ and $(\cdot|\cdot)$ are the norm and scalar product in $L^2(\Omega)$, respectively.}  \\
(3) \textit{The solution $(u,v,\partial_t v)$ satisfies the estimates}
\begin{align*}
   &\|u(t)\|_{H^1}^2 + \|v(t)\|_2^2 + \|\partial_t v(t)\|_2^2 
   + \|(-\Delta)^{-1/2} \partial_t v(t)\|_2^2 \leq C' \exp (C_0^2 \|\varphi\|_2^2 |t|),  \tag{1.1}  \\
   &\|u(t)\|_{H^2}^2 + \|v(t)\|_{H^1}^2 + \|\partial_t v(t)\|_{H^1}^2 
   \leq C'' \exp \left( C'' \exp \left( \frac{1}{2} C_0^2 \|\varphi\|_2^2 |t| \right) \right),  \tag{1.2}   
\end{align*}
\textit{for any $t \in \mathbb R$, where $C'$ depends on $\|\varphi\|_{H^1}, \|\psi_0\|_2, \|\psi_1\|_2, \|(-\Delta)^{-1/2} \psi_1\|_2$, $C''$ depends on $\|\varphi\|_{H^2}, \|\psi_0\|_{H^1}, \|\psi_1\|_{H^1}$, and $C_0$ is the best constant in the Gagliardo-Nirenberg inequality}
\begin{align*}
   \|u\|_4^2 \leq C_0 \|u\|_2 \|\nabla u\|_2   \tag{1.3}
\end{align*}
\textit{for all $u \in H_0^1(\Omega)$}.  \\

\noindent
\textbf{Remark 1.}\ \ (1) The space $H_0^1(\Omega) \oplus L^2(\Omega) \oplus (L^2 \cap (-\Delta)^{1/2}L^2)(\Omega)$ is the natural energy space for (S-iB) in the sense that $E(t)$ makes sense as a real number by (1.3).   \\
(2) No smallness condition is required on $\|\varphi\|_2$.
This is a striking difference in view of the coresponding cases of the Zakharov system and cubic nonlinear Schr\"{o}dinger equation \cite{AA2,C,GGKZ,OT1,OT2,W}.   \\
(3) The $H^1 \oplus L^2 \oplus L^2$-bound (1.1) is an improvement of an estimate in \cite{OT2}, where the RHS on (1.1) is replaced by $C'(1+t^2)\exp(2C_0^2 \|\varphi\|_2^2 t^2)$.
The $H^1$-bound of exponential type (1.1) is a natural consequence of a simple Gronwall type argument, which is much simpler than that of \cite{OT2}.   \\
(4) The $H^2 \oplus H^1 \oplus H^1$-bound (1.2) is new.
The $H^2 \oplus H^1 \oplus H^1$-bound of double exponential type (1.2) is a natural consequence of a simple Gronwall type argument on a modified energy for $H^2 \oplus H^1 \oplus H^1$ \cite{CO1,G,HnO,OV,PTV} with coefficient of exponential contribution by (1.1).
The $H^2 \oplus H^1 \oplus H^1$-bound of double exponential type (1.2) arises also in the Brezis-Gallouet argument \cite{AA2,BG,OV}, whereas our proof is independent of the Brezis-Gallouet inequality.
See the proof of Theorem 1 below.    \\

We prove Theorem 1 in Section 3.
The idea of the proof is similar to that of \cite{FMO,Hm,HmO,OT1} in the sense that we prove that solutions of some regularized system by the Yosida approximation form a bounded sequence in $H^2 \oplus H^1 \oplus H^1$ and a Cauchy sequence in $H^1 \oplus L^2 \oplus L^2$.
For the $H^1 \oplus L^2 \oplus L^2$-control (1.1) by the energy without smallness of $\|\varphi\|_2$, we employ the Gronwall argument on $v$ to avoid a direct use of the sharp Gagliardo-Nirenberg inequality (1.3).
For the $H^2 \oplus H^1 \oplus H^1$-control (1.2), we introduce a modified energy for $H^2 \oplus H^1 \oplus H^1$ which closes another Gronwall argument on $(\partial_t u, \partial_t v)$.
Those two strategies require a specific technique, which is totally different from that of \cite{OT1}.

The second purpose of this paper is to study the vanishing ''improvement'' limit problem for (S-iB) as the coefficient of the linear term of the highest order tends to zero.
For that purpose we introduce $(\textrm{S-iB})_{\varepsilon}$ with $\varepsilon \in (0,1)$ as follows:
\begin{align*}
   (\text{S-iB})_{\varepsilon}
   \begin{cases}
   i \partial_t u + \Delta u = vu,    \\
   \partial_t^2 v - \Delta v - \varepsilon \Delta \partial_t^2 v = \Delta |u|^2   
   \end{cases}
\end{align*}
with the same initial data.
The proof of Theorem 1 with minor modifications implies the global existence of solutions to the Cauchy problem for $(\textrm{S-iB})_{\varepsilon}$ in the same class as in Theorem 1.
We denote by $(u_{\varepsilon}, v_{\varepsilon}, \partial_t v_{\varepsilon})$ the solutions to $(\textrm{S-iB})_{\varepsilon}$.
We prove that, under the same smallness condition as in (Z), solutions $(u_{\varepsilon}, v_{\varepsilon}, \partial_t v_{\varepsilon})$ of $(\textrm{S-iB})_{\varepsilon}$ tend to the solution $(u, v, \partial_t v)$ of (Z) as $\varepsilon \downarrow 0$.
Specifically, we prove:   \\

\noindent
\textbf{Theorem 2.}\ \ \textit{Let $(\varphi, \psi_0, \psi_1) \in D(\Delta) \oplus H_0^1(\Omega) \oplus (H_0^1 \cap (-\Delta)^{1/2}L^2)(\Omega)$ and let $\varphi$ satisfy $\|\varphi\|_2 < \sqrt{2}/C_0$, where $C_0$ is the best constant for (1.3).
Then}:   \\
\noindent
(1) \textit{The solutions $(u_{\varepsilon}, v_{\varepsilon}, \partial_t v_{\varepsilon})$ of} $(\textrm{S-iB})_{\varepsilon}$ \textit{with $\varepsilon \in (0,1)$ satisfy the estimates}
\begin{align*}
   &\|u_{\varepsilon}(t)\|_2 = \|\varphi\|_2,   \tag{1.4}   \\
   &\|\nabla u_{\varepsilon}(t)\|_2^2 + \frac{1}{2} ( \|v_{\varepsilon}(t)\|_2^2 
   + \|(-\Delta)^{-1/2}\partial_t v_{\varepsilon}(t)\|_2^2 
   + \varepsilon \|\partial_t v_{\varepsilon}(t)\|_2^2 )   \\[-3pt]
   &\leq \left( 1- \frac{C_0}{\sqrt{2}} \|\varphi\|_2 \right)^{-1} ( \|\nabla \varphi\|_2^2 
   + \frac{1}{2} (\|\psi_0\|_2^2 + \|(-\Delta)^{-1/2} \psi_1\|_2^2 + \|\psi_1\|_2^2 )    \\[-3pt]
   &\hspace{98pt} + C_0 \|\varphi\|_2 \|\nabla \varphi\|_2 \|\psi_0\|_2 ),   \tag{1.5}   \\
   &\sup_{\varepsilon \in (0,1)} (\|u_{\varepsilon}(t)\|_{H^2}^2 
   + \|v_{\varepsilon}(t)\|_{H^1}^2+ \|\partial_t v_{\varepsilon}(t)\|_2^2) 
   \leq C'' \exp (C'|t|),   \tag{1.6}
\end{align*}
\textit{where $C'$ depends on $\|\varphi\|_{H^1}, \|\psi_0\|_2, \|\psi_1\|_2, \|(-\Delta)^{-1/2} \psi_1\|_2$ and $C''$ depends on $\|\varphi\|_{H^2}, \|\psi_0\|_{H^1}, \|\psi_1\|_{H^1}$.}   \\
(2) \textit{The solutions $(u_{\varepsilon}, v_{\varepsilon}, \partial_t v_{\varepsilon})$ of} $(\textrm{S-iB})_{\varepsilon}$ \textit{with $\varepsilon \in (0,1)$ converge to the unique global strong solution $(u, v, \partial_t v)$ of} (Z) \textit{on compact intervals with values in $H^1 \oplus L^2 \oplus (-\Delta)^{1/2}L^2$, namely, for any $T>0$,}
\begin{align*}
   &\sup_{t \in [-T,T]} (\|u_{\varepsilon}(t)-u(t)\|_{H^1} 
   + \|v_{\varepsilon}(t)-v(t)\|_2    \\[-5pt]
   &\hspace{35pt} + \|(-\Delta)^{-1/2}(\partial_t v_{\varepsilon}(t)-\partial_t v(t))\|_2 ) \to 0   \tag{1.7}
\end{align*}
\textit{as $\varepsilon \downarrow 0$.}   \\

\noindent
\textbf{Remark 2.}\ \ (1) The smallness assumption $\|\varphi\|_2<\sqrt{2}/C_0$ is necessary to ensure the global existence of solutions to (Z) as well as the uniform estimate (1.5) with respect to $\varepsilon \in (0,1)$ and $t \in \mathbb R$.   \\
(2) A similar $H^2 \oplus H^1 \oplus L^2$-bound for solutions to (Z) is proved in\cite{OT1}, where $C''$ depends on $\|\varphi\|_{H^2}, \|\psi_0\|_{H^1}, \|\psi_1\|_2$.   \\

We prove Theorem 2 in Section 4.     
   

\section{Preliminaries.}

\ \ \ In this section, we collect basic estimates to be used in the proof of the main theorems below.
We use the following Gagliardo-Nirenberg inequalities
\begin{align*}
   &\|u\|_4^2 \leq C_0 \|u\|_2 \|\nabla u\|_2,   \tag{2.1}   \\
   &\|u\|_{\infty}^2 \leq C_1 \|u\|_2 \|u\|_{H^2},   \tag{2.2}
\end{align*}
the standard elliptic estimate
\begin{align*}
   \|u\|_{H^2} \leq C_2 (\|\Delta u\|_2 + \|u\|_2),   \tag{2.3}
\end{align*}
and the elementary equalities
\begin{align*}
   &\|(-\Delta)^{1/2} u\|_2 = \|\nabla u\|_2,   \tag{2.4}   \\
   &\|(1-\Delta)^{1/2} (-\Delta)^{-1/2} u\|_2^2 
   = \|u\|_2^2 + \|(-\Delta)^{-1/2} u\|_2^2.   \tag{2.5}
\end{align*}
By (2.1), we have 
\begin{align*}
   \|\nabla u\|_4^2 \leq C_0 \|\nabla u\|_2 \|u\|_{H^2}   \tag{2.6}
\end{align*}
with the same constant $C_0$.
We also use the following lemma.   \\

\noindent
\textbf{Lemma 1}(\cite{Og,OO,O})\textbf{.}\ \ \textit{There exists a constant $C$ such that the estimate}
\begin{align*}
   \|u\|_p \leq C p^{1/2} \|u\|_2^{2/p} \|\nabla u\|_2^{1-2/p}   \tag{2.7}
\end{align*}
\textit{holds for any $u\in H^1(\Omega)$ and any $p \in [2, +\infty)$.}   \\

We summarize basic properties of the Yosida approximation of the identity.   \\

\noindent
\textbf{Lemma 2}(\cite{C,OT1})\textbf{.}\ \ \textit{For any positive integer $n$, let $J_n = (I- \frac{1}{n} \Delta)^{-1}$. Then} : \\
\noindent
(1) \textit{For any positive integer $n$, $J_n$ is a bounded self-adjoint operator in $L^2(\Omega)$ with image given by}
\begin{align*}
   J_n (L^2(\Omega)) = D(\Delta) = (H^2 \cap H_0^1)(\Omega).   
\end{align*}
(2) \textit{For any $u\in L^2(\Omega), J_n u \to u$ in $L^2(\Omega)$ as $n \to \infty$.}   \\
(3) \textit{The following estimates}
\begin{align*}
   &\|J_n u\|_2 \leq \|u\|_2,   \tag{2.8}   \\[1pt]
   &\|\nabla J_n u\|_2 \leq n^{1/2} \|u\|_2,  \tag{2.9}   \\[3pt]
   &\|\nabla J_n u\|_2 \leq \|\nabla u\|_2,   \tag{2.10}   \\[3pt]
   &\|\Delta J_n u\|_2 \leq \|\Delta u\|_2.   \tag{2.11}   
\end{align*}
\textit{hold for any $n \in \mathbb Z_{>0}$.}


\section{Proof of Theorem 1.}

\ \ \ In this section, we prove Theorem 1.
We introduce the following regularized system for (S-iB):
\begin{align*}
(\text{S-iB})_{n}
   \begin{cases}
   i \partial_t u + \Delta u = J_n (J_n v \cdot J_n u),    \\
   \partial_t^2 v - \Delta v - \Delta \partial_t^2 v = J_n \Delta |J_n u|^2   
   \end{cases}
\end{align*}
with regularized initial data
\begin{align*}
   (u(0), v(0), \partial_t v(0)) = (J_n \varphi, J_n \psi_0, J_n \psi_1).
\end{align*}
Then $(\textrm{S-iB})_{n}$ are converted into the integral equations
\begin{align*}
(\text{S-iB})_{n}'
   \begin{cases}
   \displaystyle u(t) = U(t) J_n \varphi - i \int_0^t U(t-t') J_n (J_n v \cdot J_n u)(t') dt',    \\
   \displaystyle v(t) = \dot{K}(t) J_n \psi_0 + K(t) J_n \psi_1 + \int_0^t K(t-t') J_n \Delta |J_n u|^2(t') dt',
   \end{cases}
\end{align*}
where $U(t) = \exp(it\Delta), K(t) = \omega^{-1}(\sin t\omega), \dot{K}(t) = \cos t\omega$ with $\omega =(-\Delta)^{1/2} (1-\Delta)^{-1/2}$.
By the standard fixed point argument \cite{CH}, $(\textrm{S-iB})_{n}'$ have unique global solutions $(u_n, v_n, \partial_t v_n) \in C^1(\mathbb R ; D(\Delta) \oplus D(\Delta) \oplus D(\Delta))$ since the perturbation on the RHS of $(\textrm{S-iB})_{n}'$ is globally Lipschitz in $D(\Delta)$.
Moreover, $(u_n, v_n, \partial_t v_n)$ solve $(\textrm{S-iB})_{n}$.

We first show the boundedness of the sequence $((u_n, v_n, \partial_t v_n) ; n \in \mathbb Z_{>0})$ of solutions to $(\textrm{S-iB})_{n}$ with values in $H^2 \oplus H^1 \oplus H^1$.
From now on we restrict our attention to the case $t>0$ for definiteness.
Omitting explicit dependence of the time variable, we have
\begin{align*}
   \frac{d}{dt} \|u_n\|_2^2 &= 2 \textrm{Im} (i \partial_t u_n | u_n)   \\[-5pt]
        &= 2 \textrm{Im} (-\Delta u_n + J_n (J_n v_n \cdot J_n u_n) | u_n)   \\
        &= 2 \textrm{Im} (J_n v_n \cdot J_n u_n | J_n u_n) = 0,   \tag{3.1}
\end{align*}
\vspace{-22pt}
\begin{align*}
   &\frac{d}{dt} (\|\nabla u_n\|_2^2 
   +\frac{1}{2} (\|v_n\|_2^2 + \|\partial_t v_n\|_2^2 
   + \|(-\Delta)^{-1/2} \partial_t v_n \|_2^2) + (J_n v_n | |J_n u_n|^2) )   \\[1pt]
        &= 2 \textrm{Re} (\partial_t \nabla u_n | \nabla u_n) 
        + ((-\Delta)^{-1}(- \Delta - \Delta \partial_t^2 + \partial_t^2)v_n | \partial_t v_n )   \\[2pt]
        &\hspace{10pt} + (J_n \partial_t v_n | |J_n u_n|^2) 
        + (J_n v_n |\ \partial_t |J_n u_n|^2)   \\[2pt]
        &= - 2 \textrm{Re} (\partial_t u_n | -i \partial_t u_n + J_n (J_n v_n \cdot J_n u_n)) 
        - (J_n |J_n u_n|^2 | \partial_t v_n)  \\[2pt]
        &\hspace{10pt} + (\partial_t v_n | J_n |J_n u_n|^2) 
        + (J_n v_n |\ \partial_t |J_n u_n|^2)   \\[2pt]
        &= - 2 \textrm{Re} (J_n \partial_t u_n | J_n v_n \cdot J_n u_n) 
        + (J_n v_n |\ \partial_t |J_n u_n|^2) = 0.   \tag{3.2}
\end{align*}
It follows from (2.8),(3.1), and (3.2) that
\begin{align*}
   &\|u_n (t)\|_2^2 = \|u_n (0)\|_2^2 
   = \|J_n \varphi \|_2^2 \leq \|\varphi\|_2^2,   \tag{3.3}   \\
   &E_n (t) = E_n (0),   \tag{3.4}
\end{align*}
where
\begin{align*}
   E_n (t) = \|\nabla u_n\|_2^2 
   + \frac{1}{2} (\|v_n\|_2^2 + \|\partial_t v_n\|_2^2
   + \|(-\Delta)^{-1/2}\partial_t v_n\|_2^2) + (J_n v_n||J_n u_n|^2).   \tag{3.5}
\end{align*}
From now on we assume $\varphi \neq 0$, otherwise solutions $u_n$ are trivial via (3.3).
By (2.1),(2.8),(2.10), and (3.3), the last term on the RHS of (3.5) is estimated by 
\begin{align*}
   |(J_n v_n | |J_n u_n |^2)| &\leq \|J_n v_n \|_2 \|J_n u_n \|_4^2   \\
        &\leq C_0 \|J_n v_n \|_2 \|J_n u_n \|_2 \|\nabla J_n u_n\|_2   \\
        &\leq C_0 \|v_n \|_2 \|u_n \|_2 \|\nabla u_n\|_2   \\
        &\leq C_0 \|v_n\|_2 \|\varphi\|_2 \|\nabla u_n\|_2   \\
        &\leq \frac{1}{2} \|\nabla u_n\|_2^2 + \frac{1}{2} C_0^2 \|\varphi\|_2^2 \|v_n\|_2^2. \tag{3.6}
\end{align*}
Similarly,
\begin{align*}
   |(J_n v_n(0) | |J_n u_n(0)|^2)| \leq C_0 \|\psi_0\|_2 \|\varphi\|_2 \|\nabla \varphi\|_2.   \tag{3.7}
\end{align*}
By (3.4),(3.5),(3.6), and (3.7), we obtain
\begin{align*}
   &\|\nabla u_n\|_2^2 + \|v_n\|_2^2 + \|\partial_t v_n\|_2^2 
   + \|(-\Delta)^{-1/2} \partial_t v_n\|_2^2   \\
   &\leq 2 E_n(t) + C_0^2 \|\varphi\|_2^2 \|v_n\|_2^2   \\
   &= 2 E_n(0) + C_0^2 \|\varphi\|_2^2 \|v_n\|_2^2   \\
   &\leq 2 \|\nabla \varphi\|_2^2 + \|\psi_0\|_2^2 + \|\psi_1\|_2^2 
   + \|(-\Delta)^{-1/2} \psi_1\|_2^2   \\
   &\hspace{10pt} + C_0 \|\psi_0\|_2 \|\varphi\|_2 \|\nabla \varphi\|_2
   + C_0^2 \|\varphi\|_2^2 \|v_n\|_2^2.   \tag{3.8}
\end{align*}
We estimate $\|v_n\|_2^2$ as
\begin{align*}
   \|v_n(t)\|_2^2 &= \|v_n(0)\|_2^2 + 2 \int_0^t (\partial_t v_n(t') | v_n(t')) dt'    \\[-5pt]
   &\leq \|\psi_0\|_2^2 + \int_0^t (\|v_n(t')\|_2^2 + \|\partial_t v_n(t')\|_2^2) dt'.   \tag{3.9}
\end{align*}
The Gronwall argument on (3.8) with (3.9) implies
\begin{align*}
   \|\nabla u_n(t)\|_2^2 + \|v_n(t)\|_2^2 + \|\partial_t v_n(t)\|_2^2 
   + \|(-\Delta)^{-1/2} \partial_t v_n(t)\|_2^2
   \leq C_3 \exp(C_0^2 \|\varphi\|_2^2 t),   \tag{3.10}
\end{align*}
where
\begin{align*}
   C_3 &= 2 \|\nabla \varphi\|_2^2 + \|\psi_0\|_2^2 + \|\psi_1\|_2^2
   + \|(-\Delta)^{-1/2} \psi_1\|_2^2    \\
   &\hspace{10pt} + C_0 \|\psi_0\|_2 \|\varphi\|_2 \|\nabla \varphi\|_2
   + C_0^2 \|\varphi\|_2^2 \|\psi_0\|_2^2.
\end{align*}
The estimate (3.10) with (3.3) shows the boundedness of solutions to $(\textrm{S-iB})_{n}$ with values in $H^1 \oplus L^2 \oplus L^2$.
For the boundedness in $H^2 \oplus H^1 \oplus H^1$, we introduce the modified energy by
\begin{align*}
   F_n(t) = \|\partial_t u_n\|_2^2 + \frac{1}{2} ( \|\nabla v_n\|_2^2 
   + \|\nabla \partial_t v_n\|_2^2 + \|\partial_t v_n\|_2^2)
\end{align*}
and calculate its time derivative as
\begin{align*}
   F_n'(t) &= 2 \textrm{Re} (\partial_t^2 u_n|\partial_t u_n) 
   + (\partial_t \nabla v_n|\nabla v_n) + (\nabla \partial_t^2 v_n|\nabla \partial_t v_n)
   + (\partial_t^2 v_n|\partial_t v_n)   \\
   &= 2 \textrm{Im} (\partial_t(\Delta u_n + J_n(J_n v_n \cdot J_n u_n))|\partial_t u_n)
   + (\partial_t v_n|-\Delta v_n - \Delta \partial_t^2 v_n + \partial_t^2 v_n)   \\
   &= 2 \textrm{Im} ((\partial_t J_n v_n) J_n u_n))|J_n \partial_t u_n)
   + (\partial_t v_n|J_n \Delta |J_n u_n|^2)     \\
   &= 2 \textrm{Re} 
   ((\partial_t J_n v_n)J_n u_n|J_n(-\Delta u_n + J_n(J_n v_n \cdot J_n u_n)))
   + (J_n \partial_t v_n|\Delta |J_n u_n|^2)   \\
   &= (J_n \partial_t v_n|\Delta |J_n u_n|^2 
   - 2 \textrm{Re} (\overline{J_n u_n} \Delta J_n u_n))   
   + 2 \textrm{Re} ((J_n \partial_t v_n) J_n u_n|J_n^2(J_n v_n \cdot J_n u_n))   \\
   &= 2 (J_n \partial_t v_n||\nabla J_n u_n|^2)
   + 2 \textrm{Re} ((J_n \partial_t v_n) J_n u_n|J_n^2(J_n v_n \cdot J_n u_n)).   \tag{3.11}
\end{align*}
By (2.2),(2.3),(2.8),(2.9),(2.10), and (3.3), we estimate
\begin{align*}
   \|\Delta u_n\|_2 &= \|-i \partial_t u_n + J_n(J_n v_n \cdot J_n u_n)\|_2   \\
   &\leq \|\partial_t u_n\|_2 + \|J_n u_n\|_{\infty} \|J_n v_n\|_2   \\
   &\leq \|\partial_t u_n\|_2 
   + C_1^{1/2} \|J_n u_n\|_2^{1/2} \|J_n u_n\|_{H^2}^{1/2} \|v_n\|_2   \\
   &\leq \|\partial_t u_n\|_2 
   + C_1^{1/2} C_2^{1/2} \|\varphi\|_2^{1/2} (\|\Delta u_n\| + \|u_n\|_2)^{1/2} \|v_n\|_2 \\
   &\leq \|\partial_t u_n\|_2 + \frac{1}{2} (\|\Delta u_n\|_2 + \|\varphi\|_2)
   + \frac{1}{2} C_1 C_2 \|\varphi_2 \|v_n\|_2^2,
\end{align*}
from which we have
\begin{align*}
   \|\Delta u_n\|_2 \leq 2 \|\partial_t u_n\|_2 + \|\varphi\|_2 
   + C_1 C_2 \|\varphi\|_2 \|v_n\|_2^2.   \tag{3.12}
\end{align*}
Similarly,
\begin{align*}
   \|\partial_t u_n\|_2 \leq 2 \|\Delta u_n\|_2 + \|\varphi\|_2 
   + C_1 C_2 \|\varphi\|_2 \|v_n\|_2^2.   \tag{3.13}
\end{align*}
By (2.2),(2.3),(2.6),(2.8),(2.10),(2.11), and (3.3), we estimate the RHS of the last equality of (3.11) as
\begin{align*}
   &2|(J_n \partial_t v_n||\nabla J_n u_n|^2)| 
   + 2 |((J_n \partial_t v_n) J_n u_n|J_n^2 (J_n v_n \cdot J_n u_n))|   \\
   &\leq 2 \|J_n \partial_t v_n\|_2 \|\nabla J_n u_n\|_4^2 
   + 2 \|J_n \partial_t v_n\|_2 \|J_n u_n\|_{\infty} \|J_n^2 (J_n v_n \cdot J_n u_n)\|_2   \\
   &\leq 2 \|\partial_t v_n\|_2 ( C_0 \|\nabla J_n u_n\|_2 \|J_n u_n\|_{H^2} 
   + \|J_n u_n\|_{\infty}^2 \|J_n v_n\|_2 )   \\
   &\leq 2 \|\partial_t v_n\|_2 ( C_0 \|\nabla u_n\|_2 \|u_n\|_{H^2}
   + C_1 \|u_n\|_2 \|u_n\|_{H^2} \|v_n\|_2 )   \\
   &\leq 2 \|\partial_t v_n\|_2 ( C_0 \|\nabla u_n\|_2 
   + C_1 \|\varphi\|_2 \|v_n\|_2 ) \|u_n\|_{H^2}   \\
   &\leq 2^{3/2} F_n(t)^{1/2} ( C_0 \|\nabla u_n\|_2 
   + C_1 \|\varphi\|_2 \|v_n\|_2 ) C_2 (\|\Delta u_n\|_2 + \|\varphi\|_2).  \tag{3.14}
\end{align*}
By (3.10),(3.11), and (3.14), we have for any $\varepsilon >0$
\begin{align*}
   \frac{d}{dt} (F_n(t) + \varepsilon)^{1/2} 
   &\leq 2^{1/2} C_2 (C_0 \|\nabla u_n\|_2 
   + C_1 \|\varphi\|_2 \|v_n\|_2 ) (\|\Delta u_n\|_2 + \|\varphi\|_2)   \\[-3pt]
   &\leq 2^{1/2} C_2 (C_0 + C_1 \|\varphi\|_2) 
   (\|\nabla u_n\|_2 + \|v_n\|_2 ) (\|\Delta u_n\|_2 + \|\varphi\|_2)   \\
   &\leq 2 C_2 (C_0 + C_1 \|\varphi\|_2) 
   C_3^{1/2} \exp \left( \frac{1}{2} C_0^2 \|\varphi\|_2^2 t \right) 
   (\|\Delta u_n\|_2 + \|\varphi\|_2).   \tag{3.15}
\end{align*}
Integrating both sides of (3.15) and letting $\varepsilon \downarrow 0$, we obtain
\begin{align*}
   (F_n(t))^{1/2} \leq (F_n(0))^{1/2} 
   + C_4 \int_0^t \exp \left( \frac{1}{2} C_0^2 \|\varphi\|_2^2 t' \right) 
   (\|\Delta u_n(t')\|_2 + \|\varphi\|_2) dt',   \tag{3.16}
\end{align*}
where
\begin{align*}
   C_4 = 2 C_2 (C_0 + C_1 \|\varphi\|_2) C_3^{1/2} + C_0^2 \|\varphi\|_2^2
\end{align*}
and we have added the last term for the Gronwall argument below.
By (3.10),(3.12), and (3.16), for any $\varphi \neq 0$ we obtain
\begin{align*}
   \|\Delta u_n(t)\|_2 &\leq 2 (F_n(t))^{1/2} + \|\varphi\|_2 
   + C_1 C_2 \|\varphi\|_2 \|v_n(t)\|_2^2   \\
   &\leq 2 (F_n(0))^{1/2} + \|\varphi\|_2 
   + C_1 C_2 C_3 \|\varphi\|_2 \exp(C_0^2 \|\varphi\|_2^2 t)   \\[-3pt]
   &\hspace{10pt} + C_4 \|\varphi\|_2 
   \int_0^t \exp \left( \frac{1}{2} C_0^2 \|\varphi\|_2^2 t' \right) dt'  \\[-3pt]
   &\hspace{10pt} + C_4 \int_0^t \exp \left( \frac{1}{2} C_0^2 \|\varphi\|_2^2 t' \right) 
   \|\Delta u(t')\|_2 dt'   \\[-3pt]
   &\leq 2 (F_n(0))^{1/2} + \|\varphi\|_2 + \left( C_1 C_2 C_3 \|\varphi\|_2 
   + \frac{C_4}{C_0^2 \|\varphi\|_2} \right) \exp(C_0^2 \|\varphi\|_2^2 t)   \\[-3pt]
   &\hspace{10pt} + C_4 \int_0^t \exp \left( \frac{1}{2} C_0^2 \|\varphi\|_2^2 t' \right) 
   \|\Delta u(t')\|_2 dt',   \tag{3.17} 
\end{align*}
where we have used the following inequality
\begin{align*}
   \int_0^t e^{\frac{1}{2} \lambda t'} dt' = \frac{2}{\lambda} (e^{\frac{1}{2} \lambda t}-1) 
   \leq \frac{1}{\lambda} e^{\lambda t}
\end{align*}
with $\lambda, t>0$.
Then the Gronwall argument on (3.17) yields
\begin{align*}
   \|\Delta u_n(t)\|_2 \leq C_5 
   \exp \left( \frac{2C_4}{C_0^2 \|\varphi\|_2^2} 
   \exp \left( \frac{1}{2} C_0^2 \|\varphi\|_2^2 t \right) \right),    \tag{3.18}
\end{align*}
where
\begin{align*}
   C_5 = 4 \sup_{n \geq 1} (F_n(0))^{1/2} + 2 (1+C_1 C_2 C_3) \|\varphi\|_2 
   + \frac{2C_4}{C_0^2 \|\varphi\|_2}
\end{align*}
and $F_n(0)$ is bounded uniformly in $n$ as
\begin{align*}
   F_n(0) &= \|\Delta J_n \varphi - J_n(J_n \psi_0 \cdot J_n \varphi)\|_2^2
   + \frac{1}{2} (\|\nabla J_n \psi_0\|_2^2 
   + \|\nabla J_n \psi_1\|_2^2 + \|J_n \psi_1\|_2^2)   \\
   &\leq (\|\Delta \varphi\|_2 + C_1^{1/2} \|\varphi\|_{H^2} \|\psi_0\|_2)^2
   + \frac{1}{2} (\|\nabla \psi_0\|_2^2 + \|\psi_1\|_{H^1}^2).
\end{align*}
By (3.16) and (3.18), we obtain the following uniform bound for $(v_n, \partial_t v_n)$ in $\dot{H}^1 \oplus H^1$
\begin{align*}
   \|\nabla v_n(t)\|_2 + \|\partial_t v_n(t)\|_2 
   \leq 2 (1+C_4) C_5 \exp \left( \frac{4C_4}{C_0^2 \|\varphi\|_2^2} 
   \exp \left( \frac{1}{2} C_0^2 \|\varphi\|_2^2 t \right) \right).   \tag{3.19}
\end{align*}
Collecting (3.10),(3.18), and (3.19), we obtain the boundedness of solutions $(u_n, v_n, \partial_t v_n)$ of $(\textrm{S-iB})_{n}$ on compact time-intervals with values in $H^2 \oplus H^1 \oplus H^1$, as required.

The proof of the convergence of $(u_n, v_n, \partial_t v_n)$ on compact intervals with values in $H^1 \oplus L^2 \oplus (L^2 \cap (-\Delta)^{1/2} L^2)$ is almost similar to that of \cite{OT1} and omitted.
Moreover, the existence and uniqueness of global solutions $(u, v, \partial_t v)$ to (S-iB) in the class
\begin{align*}
   &u \in (L_{loc}^{\infty} \cap C_w) (\mathbb R ; D(\Delta)) 
   \cap C(\mathbb R ; H_0^1(\Omega)) \cap C^1(\mathbb R ; H^{-1}(\Omega)),   \\
   &\partial_t u \in (L_{loc}^{\infty} \cap C_w) (\mathbb R ; L^2(\Omega) )   \\
   &v \in C(\mathbb R ; H_0^1(\Omega)) \cap C^1(\mathbb R ; L^2(\Omega)) 
   \cap C^2(\mathbb R ; H^{-1}(\Omega)),   \\
   &(-\Delta)^{-1/2} \partial_t v \in C(\mathbb R ; L^2(\Omega)),
\end{align*}
and the convergence of $(u_n, v_n, \partial_t v_n)$ to $(u, v, \partial_t v)$ on compact time-intervals with values in $H^1 \oplus L^2 \oplus (L^2 \cap (-\Delta)^{1/2} L^2)$ with weak convergence in $H^2 \oplus H^1 \oplus H^1$ follow in the same way as in the proof of Theorem 1 of \cite{OT1}.
Therefore, Part (2) of Theorem 1 follows from (3.3) and (3.4) and Part (3) follows from (3.3),(3.10),(3.18), and (3.19).
It remains to prove that $\partial_t v \in C^1(\mathbb R ; H_0^1(\Omega))$.
Let $t_0 \in \mathbb R$ and let $I=[a,b]$ an interval with $t_0 \in I$.
By (S-iB), we write
\begin{align*}
   \partial_t v(t) - \partial_t v(t_0) = \int_{t_0}^t \partial_t^2 v(t') dt'
   = \int_{t_0}^t (-\Delta) (1-\Delta)^{-1} (v+|u|^2)(t') dt',\ t \in I
\end{align*}
and estimate both sides in $H^1$ as
\begin{align*}
   \|\partial_t v(t) - \partial_t v(t_0)\|_{H^1} 
   &\leq \left| \int_{t_0}^t \|(v+|u|^2)(t')\|_{H^1} dt' \right|   \\
   &\leq \left| \int_{t_0}^t (\|v(t')\|_{H^1} + \|u(t')\|_4^2 
   + 2 \|u(t')\|_4 \|\nabla u(t')\|_4 dt' \right|   \\
   &\leq |t-t_0| \sup_{t' \in I} ( \|v(t')\|_{H^1} 
   + C_0 \|\nabla u(t')\|_2 (\|\varphi\|_2 + \|u(t')\|_{H^2}) ),
\end{align*}
which tends to zero as $t \to t_0$.
This completes the proof of Theorem 1.


\section{Proof of Theorem 2.} 

\ \ \ In this section, we prove Theorem 2.
Let $(\varphi, \psi_0, \psi_1) \in D(\Delta) \oplus H_0^1(\Omega) \oplus (H_0^1 \cap (-\Delta)^{1/2} L^2)(\Omega)$ and let $\|\varphi\|_2<\sqrt{2}/C_0$.
The proof of Theorem 1 with minor modification by replacing $\omega = (-\Delta)^{1/2}(1-\Delta)^{-1/2}$ by $\omega_{\varepsilon}= (-\Delta)^{1/2}(1-\varepsilon \Delta)^{-1/2}$ in $(\textrm{S-iB})_{n}'$ implies the existence and uniqueness of global solutions $(u_{\varepsilon}, v_{\varepsilon}, \partial_t v_{\varepsilon})$ to the Cauchy problem for $(\textrm{S-iB})_{\varepsilon}$ with $(u_{\varepsilon}(0), v_{\varepsilon}(0), \partial_t v_{\varepsilon}(0)) = (\varphi, \psi_0, \psi_1)$ in the same class as in Theorem 1.
We denote by $(u, v, \partial_t v)$ the global solution to (Z) with the same initial data $(\varphi, \psi_0, \psi_1)$(see \cite{OT1}). 
In the same way as in (3.3) and (3.4), we have conservation laws of charge and energy for $(\textrm{S-iB})_{\varepsilon}$ as
\begin{align*}
   &\|u_{\varepsilon}(t)\|_2^2 = \|\varphi\|_2^2,   \tag{4.1}   \\
   &E_{\varepsilon}(t) = E_{\varepsilon}(0),   \tag{4.2}
\end{align*}
where
\begin{align*}
   E_{\varepsilon}(t) = \|\nabla u_{\varepsilon}\|_2^2 + \frac{1}{2} ( \|v_{\varepsilon}\|_2^2 
   + \|(-\Delta)^{-1/2} \partial_t v_{\varepsilon}\|_2^2 
   + \varepsilon \|\partial_t v_{\varepsilon}\|_2^2 )
   + (v_{\varepsilon}||u_{\varepsilon}|^2).   \tag{4.3}
\end{align*}
We estimate the last term on the RHS of (4.3) as
\begin{align*}
   |(v_{\varepsilon}||u_{\varepsilon}|^2)| 
   &\leq \|v_{\varepsilon}\|_2 \|u_{\varepsilon}\|_4^2   \\
   &\leq C_0 \|v_{\varepsilon}\|_2 \|u_{\varepsilon}\|_2 \|\nabla u_{\varepsilon}\|_2   \\
   &= C_0 \|\varphi\|_2 \|v_{\varepsilon}\|_2 \|\nabla u_{\varepsilon}\|_2   \\
   &\leq \frac{C_0}{\sqrt{2}} \|\varphi\|_2 \left( \|\nabla u_{\varepsilon}\|_2^2 
   + \frac{1}{2} \|v_{\varepsilon}\|_2^2 \right),
\end{align*}
which yields the $\dot{H}^1 \oplus L^2 \oplus (L^2 \cap (-\Delta)^{1/2}L^2)$-bound
\begin{align*}
   &\|\nabla u_{\varepsilon}\|_2^2 + \frac{1}{2} ( \|v_{\varepsilon}\|_2^2 
   + \|(-\Delta)^{-1/2} \partial_t v_{\varepsilon}\|_2^2 
   + \varepsilon \|\partial_t v_{\varepsilon}\|_2^2 )   \\[-3pt]
   &\leq \left( 1-\frac{C_0}{\sqrt{2}} \|\varphi\|_2 \right)^{-1} E_{\varepsilon}(0) 
   \leq C_6,   \tag{4.4}
\end{align*}
where
\begin{align*}
   C_6 &= \left( 1-\frac{C_0}{\sqrt{2}} \|\varphi\|_2 \right)^{-1} (\|\nabla \varphi\|_2^2 
   + \frac{1}{2} (\|\psi_0\|_2^2 + \|\psi_1\|_2^2 + \|(-\Delta)^{-1/2} \psi_1\|_2^2)  \\[-5pt]
   &\hspace{100pt} + C_0 \|\varphi\|_2 \|\nabla \varphi\|_2 \|\psi_0\|_2 ).
\end{align*}
For the boundedness of $(u_{\varepsilon}, v_{\varepsilon}, \partial_t v_{\varepsilon})$ in $H^2 \oplus H^1 \oplus H^1$, we introduce the modified energy by
\begin{align*}
   F_{\varepsilon}(t) = \|\partial_t u_{\varepsilon}\|_2^2 
   + \frac{1}{2} (\|\nabla v_{\varepsilon}\|_2^2 + \|\partial_t v_{\varepsilon}\|_2^2 
   + \varepsilon \|\nabla \partial_t v_{\varepsilon}\|_2^2) + \|\varphi\|_2^2.
\end{align*}
In the same way as in (3.11), its time derivative is calculated as
\begin{align*}
   F_{\varepsilon}'(t) = 2 (\partial_t v_{\varepsilon} | |\nabla u_{\varepsilon}|^2 
   + |u_{\varepsilon}|^2 v_{\varepsilon})   \tag{4.5}
\end{align*}
and estimated as
\begin{align*}
   |F_{\varepsilon}'(t)| &\leq 2 \|\partial_t v_{\varepsilon}\|_2 
   ( C_0 \|\nabla u_{\varepsilon}\|_2 + C_1 \|\varphi\|_2 \|v_{\varepsilon}\|_2 )
   \|u_{\varepsilon}\|_{H^2}    \\
   &\leq 4 F_{\varepsilon}(t)^{1/2} (C_0 + C_1 \|\varphi\|_2) C_2 C_6^{1/2} 
   (\|\Delta u_{\varepsilon}\|_2 + \|\varphi\|_2 ).   \tag{4.6}
\end{align*}
As in (3.12), the last factor is estimated by
\begin{align*}
   \|\Delta u_{\varepsilon}\|_2 + \|\varphi\|_2 
   &\leq 2 \|\partial_t u_{\varepsilon}\|_2 + 2 \|\varphi\|_2 
   + C_1 C_2 \|\varphi\|_2 \|v_{\varepsilon}\|_2^2   \\
   &\leq 2 \|\partial_t u_{\varepsilon}\|_2 + 2 (1+C_1 C_2 C_6) \|\varphi\|_2   \\
   &\leq 2^{3/2} (1+C_1 C_2 C_6) F_{\varepsilon}(t)^{1/2}.   \tag{4.7}
\end{align*}
By (4.6) and (4.7), we obtain
\begin{align*}
   F_{\varepsilon}(t) \leq F_{\varepsilon}(0) \exp(C_7 t),   \tag{4.8}
\end{align*}
where
\begin{align*}
   C_7 = 2^{5/2} (C_0 + C_1 \|\varphi\|_2) C_2 C_6^{1/2} (1+ C_1 C_2 C_6),
\end{align*}
which together with (4.7) yields the $H^2 \oplus H^1 \oplus H^1$-bound
\begin{align*}
   \|\Delta u_{\varepsilon}\|_2^2 + \|\nabla v_{\varepsilon}\|_2^2 
   + \|\partial_t v_{\varepsilon}\|_2^2 
   + \varepsilon \|\nabla \partial_t v_{\varepsilon}\|_2^2
   \leq C_8 \exp(C_7 t),   \tag{4.9}
\end{align*}
where
\begin{align*}
   C_8 = 16 (1+C_1 C_2 C_6)^2 (4 \|\Delta \varphi\|_2^2 
   + (1+C_1^2C_2^2) \|\varphi\|_2^2 + \|\nabla \psi_0\|_2^2 
   + \|\psi_1\|_2^2 + \|\nabla \psi_1\|_2^2).
\end{align*}
This completes the proof of Part (1).

We turn to the proof of Paart (2).
By $(\textrm{S-iB})_{\varepsilon}$ and (Z), the differences of coresponding solutions are written in the form
\begin{align*}
   &i \partial_t (u_{\varepsilon}-u) + \Delta (u_{\varepsilon}-u) 
   = v_{\varepsilon} u_{\varepsilon} -vu
   = (v_{\varepsilon}-v)u_{\varepsilon} + (u_{\varepsilon}-u)v,   \tag{4.10}   \\
   &\partial_t^2 (v_{\varepsilon}-v) - \Delta (v_{\varepsilon}-v) 
   - \varepsilon \Delta \partial_t^2 v_{\varepsilon} 
   = \Delta (|u_{\varepsilon}|^2-|u|^2).   \tag{4.11}
\end{align*}
By (4.10), we have
\begin{align*}
   \frac{d}{dt} \|u_{\varepsilon}-u\|_2^2 
   &= 2 \textrm{Im} ((v_{\varepsilon}-v)u|u_{\varepsilon}-u)   \\[-5pt]
   &\leq 2 \|u\|_{\infty} \|v_{\varepsilon}-v\|_2 \|u_{\varepsilon}-u\|_2   \\
   &\leq 2 \|u\|_{H^2} \|v_{\varepsilon}-v\|_2 \|u_{\varepsilon}-u\|_2   \\
   &\leq C'' \exp(C't) \|v_{\varepsilon}-v\|_2 \|u_{\varepsilon}-u\|_2.   \tag{4.12}
\end{align*}
Here and hereafter, $C'$ and $C''$ are as in Part (1) and written by the same symbol to denote different constants independent of $t$ and $\varepsilon$.
By (4.10) and (4.11), we have
\begin{align*}
   &\frac{d}{dt} ( \|\nabla (u_{\varepsilon}-u)\|_2^2 + \frac{1}{2} (\|v_{\varepsilon}-v\|_2^2 
   + \|(-\Delta)^{-1/2} \partial_t (v_{\varepsilon}-v)\|_2^2 
   + \varepsilon \|\partial_t (v_{\varepsilon}-v)\|_2^2) )   \\
   &= 2 \textrm{Re} (\nabla \partial_t (u_{\varepsilon}-u) | \nabla (u_{\varepsilon}-u))
   + (\partial_t (v_{\varepsilon}-v) | v_{\varepsilon}-v)   \\
   &\hspace{11pt} + (\partial_t (v_{\varepsilon}-v) | (-\Delta)^{-1} 
   (\partial_t^2 -\varepsilon \Delta \partial_t^2) (v_{\varepsilon}-v))   \\
   &= -2 \textrm{Re} (\partial_t (u_{\varepsilon}-u) | \Delta (u_{\varepsilon}-u)) 
   + (\partial_t (v_{\varepsilon}-v) | (-\Delta)^{-1} 
   (\partial_t^2 - \Delta - \varepsilon \Delta \partial_t^2) (v_{\varepsilon}-v))   \\
   &= -2 \textrm{Re} (\partial_t (u_{\varepsilon}-u) | 
   (v_{\varepsilon}-v)u_{\varepsilon} + (u_{\varepsilon}-u)v) 
   - (\partial_t (v_{\varepsilon}-v) | \|u_{\varepsilon}|^2-|u|^2)   \\
   &\hspace{11pt} - \varepsilon (\partial_t (v_{\varepsilon}-v) | \partial_t^2 v )   \\
   &= -2 \textrm{Re} (\partial_t (u_{\varepsilon}-u) | (v_{\varepsilon}-v)u_{\varepsilon}) 
   - (\partial_t |u_{\varepsilon}-u|^2|v)   \\[-5pt]
   &\hspace{11pt} - (\partial_t (v_{\varepsilon}-v) | 
   2 \textrm{Re} (u_{\varepsilon}-u) \overline{u_{\varepsilon}} - |u_{\varepsilon}-u|^2)
   - \varepsilon (\partial_t v_{\varepsilon} | \partial_t^2 v)
   + \frac{\varepsilon}{2} \frac{d}{dt} \|\partial_t v\|_2^2   \\[-3pt]
   &= - \frac{d}{dt} ( 2 \textrm{Re} (u_{\varepsilon}-u | (v_{\varepsilon}-v)u_{\varepsilon}) 
   - (|u_{\varepsilon}-u|^2|v) )  
   + 2 \textrm{Re} (u_{\varepsilon}-u | 
   (v_{\varepsilon}-v) \partial_t u_{\varepsilon})  \\[-5pt]
   &\hspace{11pt} - (|u_{\varepsilon}-u|^2 | \partial_t v )
   - \varepsilon (\partial_t v_{\varepsilon} | \Delta (v + |u|^2))
   + \frac{\varepsilon}{2} \frac{d}{dt} \|\partial_t v\|_2^2,
\end{align*}
from which we have
\begin{align*}
   D_{\varepsilon}'(t) 
   = 2 \textrm{Re} (u_{\varepsilon}-u | (v_{\varepsilon}-v) \partial_t u_{\varepsilon}) 
   - (|u_{\varepsilon}-u|^2 | \partial_t v) 
   + \varepsilon (\nabla \partial_t v_{\varepsilon} | \nabla (v + |u|^2)),   \tag{4.13}
\end{align*}
where
\begin{align*}
   D_{\varepsilon}(t) &:= \|\nabla (u_{\varepsilon}-u)\|_2^2 
   + \frac{1}{2} (\|v_{\varepsilon}-v\|_2^2 
   + \|(-\Delta)^{-1/2} \partial_t (v_{\varepsilon}-v)\|_2^2 
   + \varepsilon \|\partial_t (v_{\varepsilon}-v)\|_2^2   \\[-3pt]
   &\hspace{11pt} - \varepsilon \|\partial_t v\|_2^2) 
   +2 \textrm{Re} (u_{\varepsilon}-u | (v_{\varepsilon}-v)u_{\varepsilon}) 
   - (|u_{\varepsilon}-u|^2|v).
\end{align*}
Indefinite terms in $D_{\varepsilon}(t)$ are estimated as
\begin{align*}
   |\textrm{Re} (u_{\varepsilon}-u | (v_{\varepsilon}-v)u_{\varepsilon})|  
   &\leq \|u_{\varepsilon}-u\|_2 \|v_{\varepsilon}-v\|_2 \|u_{\varepsilon}\|_{\infty}   \\
   &\leq C'' \|u_{\varepsilon}-u\|_2 \|v_{\varepsilon}-v\|_2,   \\
   |(|u_{\varepsilon}-u|^2|v)|
   &\leq C_0 \|u_{\varepsilon}-u\|_2 \|\nabla (u_{\varepsilon}-u)\|_2 \|v\|_2   \\
   &\leq C' \|u_{\varepsilon}-u\|_2 \|\nabla (u_{\varepsilon}-u)\|_2,
\end{align*}
which in turn yield
\begin{align*}
   D_{\varepsilon}(t) 
   &\leq 2 \|\nabla (u_{\varepsilon}-u)\|_2^2 
   + \|v_{\varepsilon}-v\|_2^2 
   + \|(-\Delta)^{-1/2} \partial_t (v_{\varepsilon}-v)\|_2^2 
   + \varepsilon \|\partial_t (v_{\varepsilon}-v)\|_2^2   \\[-3pt]
   &\hspace{11pt} - \varepsilon \|\partial_t v\|_2^2
   + C'' \|u_{\varepsilon}-u\|_2^2   \\
   &\leq 4 D_{\varepsilon}(t) + C'' \|u_{\varepsilon}-u\|_2^2.   \tag{4.14}
\end{align*}
We estimate contribution of each term on the RHS of (4.13), separately.
The first term is estimated by
\begin{align*}
   &|2 \textrm{Re} (u_{\varepsilon}-u | (v_{\varepsilon}-v) \partial_t u_{\varepsilon})| 
   \leq 2 \|u_{\varepsilon}-u\|_{2/\delta} \|v_{\varepsilon}-v\|_2^{1-2\delta}
   \|v_{\varepsilon}-v\|_4^{2\delta} \|\partial_t u\|_2   \\
   &\leq 2 C_0^{\delta} \|u_{\varepsilon}-u\|_{2/\delta} \|v_{\varepsilon}-v\|_2^{1-\delta}
   \|\nabla (v_{\varepsilon}-v)\|_2^{\delta} \|\partial_t u\|_2   \\
   &\leq C \delta^{-1/2} \|u_{\varepsilon}-u\|_2^{\delta} 
   \|\nabla (u_{\varepsilon}-u)\|_2^{1-\delta} 
   \|v_{\varepsilon}-v\|_2^{1-\delta} 
   (\|\nabla v_{\varepsilon}\|_2 + \|\nabla v\|_2)^{\delta} \|\partial_t u\|_2   \\
   &\leq C'' \delta^{-1/2} \|u_{\varepsilon}-u\|_{H^1} 
   \|v_{\varepsilon}-v\|_2^{1-\delta} \exp(C't),   \tag{4.15}
\end{align*}
where we have used the H\"{o}lder inequality with $1 = \frac{\delta}{2} + \frac{1-2\delta}{2} + \frac{2\delta}{4} + \frac{1}{2}$ and $0<\delta<\frac{1}{2}$, (2.7) with $p=\frac{2}{\delta}$, and (4.9), and $C'$ and $C''$ are independent of $\varepsilon$ and $\delta$.
The second and third terms are estimated respectively by
\begin{align*}
   |(|u_{\varepsilon}-u|^2 | \partial_t v)| 
   &\leq \|u_{\varepsilon}-u\|_4^2 \|\partial_t v\|_2   \\
   &\leq C'' \|u_{\varepsilon}-u\|_{H^1}^2 \exp(C't),   \tag{4.16}   \\
   |\varepsilon (\nabla \partial_t v_{\varepsilon} | \nabla (v + |u|^2))| 
   &\leq \varepsilon \|\nabla \partial_t v_{\varepsilon}\|_2 
   (\|\nabla v\|_2 + \|\nabla |u|^2\|_2)   \\
   &\leq C'' \varepsilon^{1/2} \exp(C't),   \tag{4.17}
\end{align*}
where we have used (4.9) and $C'$ and $C''$ are independent of $\varepsilon$ and $\delta$.
Collecting (4.12)-(4.17), we obtain
\begin{align*}
   &\frac{d}{dt} (D_{\varepsilon}(t) + \|u_{\varepsilon}-u\|_2^2 + \varepsilon^{1/2})   \\
   &\leq C'' \exp(C't) ( \delta^{-1/2} (D_{\varepsilon}(t) 
   + \|u_{\varepsilon}-u\|_2^2)^{1-\delta/2} + (D_{\varepsilon}(t) 
   + \|u_{\varepsilon}-u\|_2^2) + \varepsilon^{1/2} ),   \tag{4.18}
\end{align*}
where $C'$ and $C''$ are independent of $\varepsilon$ and $\delta$.
By (4.8), we have
\begin{align*}
   &\frac{d}{dt} (D_{\varepsilon}(t) + \|u_{\varepsilon}-u\|_2^2 
   + \varepsilon^{1/2})^{\delta/2}   \\[-5pt]
   &= \frac{\delta}{2} (D_{\varepsilon}(t) + \|u_{\varepsilon}-u\|_2^2 
   + \varepsilon^{1/2})^{\delta/2-1}
   \frac{d}{dt} (D_{\varepsilon}(t) + \|u_{\varepsilon}-u\|_2^2 + \varepsilon^{1/2})   \\[-1pt]
   &\leq C'' \exp(C't) (\delta^{1/2} 
   + \delta (D_{\varepsilon}(t) + \|u_{\varepsilon}-u\|_2^2 
   + \varepsilon^{1/2})^{\delta/2}).   \tag{4.19}
\end{align*}
Integrating both sides of the differential inequality (4.19) in $t$ and applying the Gronwall argument to the resulting inequality, we obtain
\begin{align*}
   D_{\varepsilon}(t) + \|u_{\varepsilon}-u\|_2^2 + \varepsilon^{1/2} 
   \leq C'' \exp(C't) ((D_{\varepsilon}(0) + \varepsilon^{1/2})^{\delta/2} 
   + \delta^{1/2} t)^{2/\delta}.   \tag{4.20}
\end{align*}
For any $\delta \in (0, \frac{1}{2})$, there exists $\varepsilon_0 >0$ such that $D_{\varepsilon}(0) + \varepsilon^{1/2} < (\frac{1}{2})^{2/\delta}$ for any $\varepsilon \in (0, \varepsilon_0]$ and therefore (4.20) implies
\begin{align*}
   D_{\varepsilon}(t) + \|u_{\varepsilon}-u\|_2^2 + \varepsilon^{1/2} 
   \leq C'' \exp(C't) \left( \frac{1}{2} + \delta^{1/2} t \right)^{2/\delta}
\end{align*}
for any $\varepsilon \in (0, \varepsilon_0]$.
Then it follows that
\begin{align*}
   \limsup_{\varepsilon \downarrow 0} \sup_{t \in [0,T]} 
   (D_{\varepsilon}(t) + \|u_{\varepsilon}(t)-u(t)\|_2^2) 
   \leq C'' \exp(C'T) \left( \frac{1}{2} + \delta^{1/2} T \right)^{2/\delta}   \tag{4.21}
\end{align*}
for any $T>0$ and $\delta \in (0,\frac{1}{2})$.
For any $T>0$, letting $\delta \downarrow 0$ in (4.21) yields
\begin{align*}
   \lim_{\varepsilon \downarrow 0} \sup_{t \in [0,T]} 
   (D_{\varepsilon}(t) + \|u_{\varepsilon}(t)-u(t)\|_2^2) = 0.   \tag{4.22}
\end{align*}
The proof of Part (2) is complete in view of (4.14) and (4.22)   \\

\noindent
\textbf{Acknowledgments.}\ \ This work is partially supported by Grant-in Aid for Scientific Research (A) 19H00644, JSPS.



\end{document}